\documentclass[leqno,12pt]{amsart}

\usepackage{amssymb}
\usepackage{amscd}              
\usepackage[all]{xy}            
\xyoption{curve}
\CompileMatrices

\newcommand{\BG}{\text{BG}}
\newcommand{\BGL}{\text{BGL}}
\newcommand{\bd}{{\bf d}}
\newcommand{\bgamma}{\text{\boldmath $\gamma$}}

\newcommand{\btau}{\text{\boldmath $\tau$}}

\newcommand{\bpi}{\text{\boldmath $\pi$}}
\newcommand{\C}{\mathbb{C}}

\newcommand{\comp}{\circ}

\newcommand{\del}{\partial}

\newcommand{\E}{{\mathcal{E}}}

\newcommand{\Ee}{\mathbb{E}}

\newcommand{\ev}{\operatorname{ev}}

\newcommand{\F}{{\mathcal{F}}}

\newcommand{\GL}{\text{GL}}

\newcommand{\GVB}{\text{GVB}}
\newcommand{\GVBD}{\text{GVBD}}

\newcommand{\id}{\text{id}}

\newcommand{\injto}{\hookrightarrow}
\newcommand{\inv}{\text{inv}}
\newcommand{\isomorph}{\cong}           
\newcommand{\isomto}{\overset{\sim}{\to}}

\newcommand{\Isom}{\text{Isom}}

\newcommand{\KGL}{\text{KGL}}

\newcommand{\Mb}{\overline{M}}

\newcommand{\N}{{\mathbb{N}}}                    

\newcommand{\Ob}{\overline{\text{\bf O}}}
\newcommand{\ol}{\overline}
\newcommand{\Oo}{{\mathcal{O}}}                  

\newcommand{\ori}{\text{or}}

\newcommand{\Pp}{{\mathbb{P}}}                   
\newcommand{\pr}{\text{pr}}

\newcommand{\st}{\text{st}}

\newcommand{\tensor}{\otimes}

\newcommand{\To}{\longrightarrow}

\newcommand{\Z}{{\mathbb{Z}}}                    

\newtheorem{theorem}{Theorem}[section]

\theoremstyle{definition}
\newtheorem{definition}[theorem]{Definition}
\newtheorem{remark}[theorem]{Remark}

\makeatletter
\def\overunderbraces #1#2#3{{%
 \baselineskip\z@skip \lineskip4\p@ \lineskiplimit4\p@
 \displaystyle  
 \setbox\z@\vbox{\ialign{&\hfil${}##{}$\hfil\cr
   \global\let\br\br@label #1\cr 
   \global\let\br\br@down #1\cr   
   #2\cr 
 }}
 \dimen@-\ht\z@ 
 \setbox\z@\vbox{\ialign{&\hfil${}##{}$\hfil\cr
   \global\let\br\br@label #1\cr 
   \global\let\br\br@down #1\cr   
   #2\cr 
   \global\let\br\br@up #3\cr 
   \global\let\br\br@label #3\cr   
 }}
 \advance\dimen@\ht\z@ 
 \lower\dimen@\hbox{\box\z@} 
}}
\def\br@up#1#2{\multispan{#1}\upbracefill}
\def\br@down#1#2{\multispan{#1}\downbracefill}
\def\br@label#1#2{\multispan{#1}\hidewidth $#2$\hidewidth}
\makeatother

\begin{document}

\title[Stable maps into the classifying space of the general linear group]
{Stable maps into the classifying space of the general linear group}
\author[Ivan Kausz]{Ivan Kausz}
\date{\today}
\address{NWF I - Mathematik, Universit\"{a}t Regensburg, 93040 Regensburg, 
Germany}
\email{ivan.kausz@mathematik.uni-regensburg.de}

\maketitle



\section{Introduction}

In this note we give a definition of stable maps into the classifying
stack $\BGL_r$ of the general linear group. 
To support our belief that the definition is the correct one, we
show that there are natural boundary morphisms between the
moduli groupoids parameterizing stable maps to $\BGL_r$.
We warn the reader beforehand that we do not as yet have applications
of our constructions. Our purpose in this note is mainly to explain
the combinatorics of stable maps into $\BGL_r$. 

Before giving an idea of what stable maps to $\BGL_r$ are, we
first recall the case of stable maps into a complex 
projective variety $V$.
Let $g\in\N_0$, and let $S$ be a finite set.
A stable $S$-pointed map of genus $g$ consists
of a prestable curve $C$ together with nonsingular points
$(x_s)_{s\in S}$ on it and a morphism $f:C\to V$ of algebraic
varieties such that the collection of these data has at most 
finitely many automorphisms.

Stable maps into $V$ are a fundamental notion
for defining Gromov-Witten invariants for $V$.
Namely, there is a Deligne-Mumford stack $\ol{M}_{g,S}(V)$ 
parameterizing all $S$-pointed stable maps of genus $g$, and
Gromov-Witten invariants of $V$ are defined by means
of intersection theory on $\ol{M}_{g,S}(V)$.

The moduli stack $\ol{M}_{g,S}(V)$ breaks up into connected
components $\ol{M}_{g,S}(V,\beta)$ where $\beta$ runs through
some monoid. Stable maps which belong to that component are
called of class $\beta$. 
Given two stable maps, say $S_i\cup\{*_i\}$-pointed, of genus $g_i$
and class $\beta_i$ for $i=1,2$ such that the points $*_i$ are 
mapped to the same point in $V$, one obtains a new stable
map which is $S_1\cup S_2$-pointed, of genus $g_1+g_2$ and
class $\beta_1+\beta_2$, by clutching together the two curves
at the points $*_i$. This yields a morphism
$$
\ol{M}_{g_1,S_1\cup\{*_1\}}(V,\beta_1)\times_V
\ol{M}_{g_2,S_2\cup\{*_2\}}(V,\beta_2)\to
\ol{M}_{g_1+g_2,S_1\cup S_2}(V,\beta_1+\beta_2)
$$
which is an instance of a boundary morphism. A general point
in its image belongs to a stable map whose underlying curve
has two smooth irreducible components joined in one double point
such that the restriction of the map to the $i$-th component
is a $S_i$-pointed stable map of genus $g_i$ and class $\beta_i$.
The combinatorial type of such a stable map is a decorated
graph containing two vertices labeled by $(g_1,\beta_1)$ 
and $(g_2,\beta_2)$  and half-edges (or ``tails'')
labeled by elements of $S_1$ and $S_2$ respectively:
$$
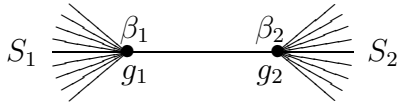
\xy
0;<1cm,0cm>:
(-1,0);(1,0) **@{-};
(-1,0)*{\bullet}; 
(1,0)*{\bullet}; 
(-1,0);(-1.766,0.642) **@{-}; 
(-1,0);(-1.766,-0.642) **@{-}; 
(1,0);(1.766,0.642) **@{-}; 
(1,0);(1.766,-0.642) **@{-}; 
(-1,0);(-1.866,0.5) **@{-}; 
(-1,0);(-1.866,-0.5) **@{-}; 
(1,0);(1.866,0.5) **@{-}; 
(1,0);(1.866,-0.5) **@{-}; 
(-1,0);(-1.939,0.342) **@{-}; 
(-1,0);(-1.939,-0.342) **@{-}; 
(1,0);(1.939,0.342) **@{-}; 
(1,0);(1.939,-0.342) **@{-}; 
(-1,0);(-1.984,0.173) **@{-}; 
(-1,0);(-1.984,-0.173) **@{-}; 
(1,0);(1.984,0.173) **@{-}; 
(1,0);(1.984,-0.173) **@{-}; 
(-1,0);(-2,-0) **@{-}; 
(1,0);(2,0) **@{-}; 
(-2.4,0)*{S_1};
(2.4,0)*{S_2};
(-0.9,-0.3)*{g_1};
(-0.9,0.3)*{\beta_1};
(0.9,-0.3)*{g_2};
(0.9,0.3)*{\beta_2};
\endxy
$$

More generally,
to each $S$-pointed stable  map of genus $g$ and class $\beta$ one
can associate some decorated graph $\tau$, its combinatorial type.
For each combinatorial type
there is a stack $\ol{M}(V,\tau)$, defined as fiber product
of some stacks $\ol{M}_{g_i,S_i}(V,\beta_i)$, 
and a boundary morphism
$$
\ol{M}(V,\tau)\to\ol{M}_{g,S}(V,\beta)
$$
generalizing the clutching described above.
The collection of boundary morphisms is responsible for
the rich algebraic structure of quantum cohomology of $V$
(cf. \cite{Manin}).

By the very definition of the stack $\BGL_r$ a morphism
from a curve $C$ to $\BGL_r$ is nothing else but a vector
bundle of rank $r$ over $C$. Thus it is clear that
the generic stable map to $\BGL_r$ should be a pair consisting
of a smooth curve together with a vector bundle of rank $r$ on it.
Now it is well known that given a family of stable curves which is 
generically smooth and has a singular special fiber, a vector
bundle defined on the complement of that fiber does not 
necessarily extend to a vector bundle over the whole family,
not even if one allows for an \'etale base change of the family.
This is just like in the situation of morphisms into a
variety $V$, where it may happen that a morphism to $V$, defined on
the complement of the special fiber, does not extend to a morphism
from the whole family to $V$, even after base change. 

However there is always a family of {\em semistable} curves 
which outside  the special fiber coincides with the given one, 
such that the
morphism extends to that new family, at least after an
\'etale base change. 
The same solution applies to vector bundles, i.e. to morphisms
into $\BGL_r$. Namely there exists a semistable family such
that the given vector bundle extends after \'etale base change.
Furthermore the new family and the extension can be chosen in such a 
way that the restriction of the vector bundle on those
chains of  projective lines in the semistable special fiber
which are contracted under stabilization
satisfies a certain condition which we call ``admissibility''
(cf. Definition \ref{adm}). Accordingly, we define $S$-pointed
stable maps of genus $g$ and degree $d$ to $\BGL_r$ as semistable
$S$-pointed curves of genus $g$
together with vector bundles of rank $r$ and degree $d$
satisfying that condition.

It is straightforward to define families of  stable $S$-pointed
maps of genus $g$ and degree $d$ to $\BGL_r$,
and pull back and isomorphisms of these. Thus we obtain a groupoid
over $\C$ which we denote by $\Mb_{g,S}(\BGL_r,d)$.

To give an idea of what the boundary morphisms are, let us first consider
the analogue of the example above. Thus we are given two curves, together
with vector bundles of rank $r$ on them, and on each of the curves
one special point. Now if we join the curves by identifying
the two special points, we need an additional datum to obtain
a vector bundle on the new curve. This datum consists in an isomorphism
between the respective fibers of two given bundles at the special points.
However isomorphisms between two vector spaces may degenerate in
families. In other words, after chosing a bases of the two vector
spaces, an isomorphism between them is given by an 
element in $\GL_r$ which is noncompact and therefore in general 
a family of isomorphisms 
cannot be extended to special points. Our solution of this problem
is a certain compactification $\KGL_r$  of $\GL_r$ which we have 
constructed in an earlier paper \cite{kgl}.

To be more precise, let us give a definition of the groupoid
$\ol{M}(\BGL_r,\btau)$, where $\btau$ is the decorated graph
$$
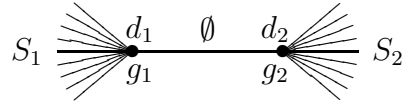
\xy
0;<1cm,0cm>:
(-1,0);(1,0) **@{-};
(-1,0)*{\bullet}; 
(1,0)*{\bullet}; 
(-1,0);(-1.766,0.642) **@{-}; 
(-1,0);(-1.766,-0.642) **@{-}; 
(1,0);(1.766,0.642) **@{-}; 
(1,0);(1.766,-0.642) **@{-}; 
(-1,0);(-1.866,0.5) **@{-}; 
(-1,0);(-1.866,-0.5) **@{-}; 
(1,0);(1.866,0.5) **@{-}; 
(1,0);(1.866,-0.5) **@{-}; 
(-1,0);(-1.939,0.342) **@{-}; 
(-1,0);(-1.939,-0.342) **@{-}; 
(1,0);(1.939,0.342) **@{-}; 
(1,0);(1.939,-0.342) **@{-}; 
(-1,0);(-1.984,0.173) **@{-}; 
(-1,0);(-1.984,-0.173) **@{-}; 
(1,0);(1.984,0.173) **@{-}; 
(1,0);(1.984,-0.173) **@{-}; 
(-1,0);(-2,-0) **@{-}; 
(1,0);(2,0) **@{-}; 
(-2.4,0)*{S_1};
(2.4,0)*{S_2};
(-0.9,-0.3)*{g_1};
(-0.9,0.3)*{d_1};
(0.9,-0.3)*{g_2};
(0.9,0.3)*{d_2};
(0,0.3)*{\emptyset};
\endxy
$$
which corresponds to the example considered above.
First we observe that one can interpret the diagonal morphism
$\Delta:\BGL_r\to\BGL_r\times\BGL_r$ as follows. There is a commutative
diagram
$$
\xymatrix{
\BGL_r \ar[rr]^(.4){\isomorph} \ar[dr]^{\Delta}
& &
\Isom(\pr_1^*\Ee,\pr_2^*\Ee) \ar[dl]
\\
& 
\BGL_r\times\BGL_r
}
$$
where the horizontal arrow is an isomorphism,
$\Ee$ is the universal vector bundle over $\BGL_r$, and 
$\pr_i:\BGL_r\times\BGL_r\to\BGL_r$ is the $i$-th projection.
We define
$$
\ol{\BGL_r}
:=\KGL(\pr_1^*\Ee,\pr_2^*\Ee)
$$
as the fiberwise compactification of the locally trivial
$\GL_r$-bundle 
$\Isom(\pr_1^*\Ee,\pr_2^*\Ee)\to\BGL_r\times\BGL_r$ and
let 
$$
\ol{\Delta}:\ol{\BGL_r}\to\BGL_r\times\BGL_r
$$
be the natural projection. With this notation the 
groupoid $\ol{M}(\BGL_r,\btau)$ is defined by the Cartesian
diagram
$$
\xymatrix{
\ol{M}(\BGL_r,\btau) \ar[r]\ar[d]
&
\ol{M}_{g_1,S_1\cup\{*_1\}}(\BGL_r,d_1)
\times
\ol{M}_{g_2,S_2\cup\{*_2\}}(\BGL_r,d_2)
\ar[d]^{\ev_1\times\ev_2}
\\
\ol{\BGL_r}\ar[r]^{\ol{\Delta}}
&
\BGL_r\times\BGL_r
}
$$
The existence of a natural boundary morphism
$$
\Mb(\BGL_r,\btau)\to\Mb_{g,S}(\BGL_r,d)
\quad.
$$
where $g=g_1+g_2$, $S=S_1\cup S_2$ and $d=d_1+d_2$ then
follows from results in \cite{degeneration} and \cite{factorization}.

Under this boundary morphism, a point of $\Mb(\BGL_r,\btau)$
which lies over the complement of $\BGL_r$ in $\ol{\BGL}_r$ is mapped
to a stable map to $\BGL_r$ whose underlying curve not stable but
only semistable, and the stabilizing map contracts a chain of
projective lines. The restriction of the underlying
bundle of the stable map to this chain of projective lines yield
a tuple of degrees and we in \S \ref{gi} below we give a careful 
description of which strata of $\ol{\BGL_r}$ (or of $\KGL_r$)
correspond to which tuple of degrees.

All the groupoids considered here are in fact Artin stacks.
Unfortunately they are neither of finite type over $\C$ nor
separated. This is to be expected since already the moduli
stack of vector bundles over a smooth curve are non-separated 
Artin stacks which are not of finite type. Still they might
be usefull e.g. for defining Gromov-Witten invariants of $\BGL_r$
if we succeed to find a kind of Shatz stratification
for them. We hope to address this question in a future paper.

It should be mentioned that A. Schmitt \cite{Sch} 
and J. Li \cite{Li} have 
constructed varieties
parameterizing  a certain subset of
all stable maps to $\BGL_r$. However they do not
define boundary morphisms  and we doubt whether
that would be possible at all in their context.

We would also like to mention work of D. Abramovich, A.Corti and 
A. Vistoli
\cite{AV,ACV} where stable maps to $\BG$ for finite
groups $G$ are defined. In \cite{ramified_bundles} we have
established some connection between their construction
an ours, but the exact relationship remains mysterious.

\section{Stable maps}

First we recall some definitions from \cite{degeneration}

\begin{definition}
A {\em chain of projective lines} 
is a connected prestable curve over a field, whose components are projective
lines and whose associated graph is linear.
An irreducible component of $R$ is called {\em extremal}, if
it intersects only one other component.
\end{definition}

\begin{definition}
\label{adm}
A vector bundle $E$ of rank $r$ on a chain $R$ of projective lines is called
{\em admissible}, if the following two properties are satisfied:
\begin{enumerate}
\item
For each irreducible component $R_i$ of $R$
there exists a number $d_i\geq 1$ such that 
the restriction of $E$ to $R_i\isomorph\Pp^1$ is isomorphic to
the bundle $d_i\Oo\oplus(r-d_i)\Oo(1)$.
\item
If a global section of $E$ vanishes in two smooth points lying in different
extremal components of $R$, then it vanishes everywhere.
\end{enumerate}
\end{definition}

It is not difficult to see that if there exists an admissible vector
bundle $E$ of rank $r$ on a chain $R$ of projective lines, then the
length of that chain is bounded by $r$ (cf. section 3 in \cite{degeneration}) .

\begin{definition}
Let $C'$ be a prestable curve over a field $k$. A 
{\em Gieseker vector bundle of rank $r$ on $C'$} is a 
pair $(C\to C', E)$, where $C$ is a prestable curve over $k$, $C\to C'$ is a 
morphism and $E$ is a rank $r$ vector bundle on $C$ 
which satisfy the following properties:
\begin{enumerate}
\item
The morphism $C\to C'$ is an isomorphism over the complement of
the singular points of $C'$.
\item
The fiber of $C\to C'$ over a singular point of $C'$ is either
a point or a chain of projective lines.
\item
The restriction of $E$ to each of the nontrivial fibers of $C\to C'$
is admissible.
\end{enumerate}
\end{definition}

Now we come to the main definition of this note:
\begin{definition}
Let $S$ be a finite set.
A {\em stable $S$-pointed map to $\BGL_r$} is a tupel
$(C,x_i\ |\ i\in S,E)$, where $(C,x_i\ |\ i\in S)$ is a prestable
$S$-pointed curve over a field and $E$ is a vector bundle on $C$ such that 
if $(C,(x_i))\to (C^{\st},(y_i))$ denotes
the stabilization morphism (cf. \cite{Manin} V.1.6) 
then $(C\to C^{\st},E)$ is a Gieseker
vector bundle of rank $r$ on $C^{\st}$.
\end{definition}

\section{Combinatorial types}
\label{ct}

Before defining the combinatorial type of a stable $S$-pointed map
to $\BGL_r$, let me first recall some definitions from \cite{Manin}.
\begin{definition}
A {\em (finite) graph} $\tau$ consists of the data 
$(F_\tau, V_\tau, \del_\tau, j_\tau)$,
where $F_\tau$ is a (finite) set (of flags), 
$V_\tau$ is a finite set (of vertices), $\del_\tau: F_\tau\to V_\tau$
is the boundary map, and $j_\tau: F_\tau\to F_\tau$ is an involution.
\end{definition}

An isomorphism $\tau\isomto\sigma$ between two graphs consists of two
bijections $F_\tau\isomto F_\sigma$, $V_\tau\isomto V_\sigma$,
compatible with $\del$ and $j$.

The fixed points
of $j_\tau$ form the set $S_\tau$ of tails,
the pairs $(f,f')$ with $f\in F_\tau$ and $f'=j_\tau(f)\neq f$ form 
the set $E^{\ori}_\tau$ of oriented edges of $\tau$.
If $e=(f,f')\in E^{\ori}$ is an oriented edge, then
$j_\tau(e):=(j_\tau(f),j_\tau(f'))$ is the edge with the
opposite orientation.
The two-element orbits of $j_\tau$ (viewed as an inversion of $F_\tau$
or of $E^{\ori}_\tau$) form the set $E$ of (un-oriented) edges of $\tau$.
For $v\in V_\tau$ let $F_\tau(v):=·\del^{-1}(v)$ denote the set
of flags starting from $v$, and let $|v|:=|F_\tau(v)|$ denote the valence of
$v$.
The topological realization of a graph $\tau$ is denoted by $||\tau||$.

\begin{definition}
A {\em modular graph} is a graph $\tau$ together with a map 
$g:V_\tau\to\N_0$, $v\mapsto g_\tau$. 
A modular graph $\tau$ is called {\em stable}, if $|v|\geq 3$ for
all $v$ with $g_v=0$ and $|v|\geq 1$ for all $v$ with $g_v=1$.
\end{definition}

There is a standard way to associate a modular graph $(\tau, g)$ 
to a prestable $S$-pointed curve $(C,x_i\ |\ i\in X)$
(cf. \cite{Manin}, III. §2).
The vertices $v$  of $\tau$ correspond to the irreducible components $C_v$ 
of $C$, and $g_v$ is the genus of $C_v$. 
The un-oriented edges $e$ of $\tau$ correspond to the singular points
$x_e$ of $C$. 
A prestable $S$-pointed curve is stable if and only if $(\tau, g)$
is stable. 

\begin{definition}
A {\em chain-type} 
is a tupel $\bd=(d_1,\dots,d_q)$ of integers $d_i\geq 1$, where
$q\geq 0$. The number $q$ is called the {\em length} of the chain-type $d$
and the integer $|\bd|:=\sum_{i=1}^qd_i$ is called its {\em degree}.
The (empty) chain-type $\bd=()$ of length zero will also be denoted by 
the symbol $\emptyset$.
\end{definition}

For a chain-type $\bd=(d_1,\dots,d_q)$ we denote by $j(\bd)$ the reversed
chain-type $(d_q,\dots,d_1)$. The reversed empty chain-type is again
empty.

\begin{definition}
A {\em chain-graph} (for rank $r$) is a tupel 
$(\tau, g, d, \bd)$ where $(\tau, g)$ is a stable modular
graph and $d$, $\bd$ are maps
\begin{eqnarray*}
d &:& V_\tau \to \Z \\
\bd &:& E_\tau^{\ori} \to \{\text{chain-types of degree $\leq r$}\}
\end{eqnarray*}
where $\bd$ has the property that
$
\bd(j_\tau(e))=j(\bd(e))
$
for every oriented edge $e$ of $\tau$.
\end{definition}

\begin{definition}
The {\em combinatorial type} of an $S$-labeled stable map
$(C, x_i\ |\ i\in S, E)$ to $\BGL_r$ is the chain-graph
$(\tau, g, d, \bd)$, where $(\tau, g)$ is the stable
modular graph associated to the stabilization $C^{\st}$ of
$C$ (in particular, $S_\tau$ is identified with $S$)
and the maps $d$, $\bd$ are characterized by the following properties. 
\begin{enumerate}
\item
Let $v\in V_\tau$ and let $C_v$ be the associated component
of the stabilized curve $C^{\st}$. There exists a unique component of
$C$ which is mapped isomorphically onto $C_v$. By a slight abuse of 
notation we denote that component again by the symbol $C_v$.
We require that the restriction of $E$ to $C_v$ is of degree $d(v)$.
\item
Let $e=(f,f')$ be an oriented edge of $\tau$ and let $x_e\in C^{\st}$ 
be the corresponding singular point. If the fiber of $C\to C^{\st}$
over $x_e$ consists of just one point, then $\bd(e)=\emptyset$. Else
the fiber $R$ is a chain of length $q\geq 1$ of projective lines. 
To the flags $f$ and $f'$ there are associated components $C_f$ 
$C_{f'}$ of  $C^{\st}$ together with points $x_f\in C_f$
and $x_{f'}\in C_{f'}$. By abuse of notation we denote by $(C_f,x_f)$
and $(C_{f'},x_{f'})$ also the pointed components of $C$ which are 
mapped isomorphically to the corresponding pointed components of 
$C^{\st}$.
Let $R_1,\dots,R_q$ be the irreducible components of the chain $R$ 
numbered in such
a way that $R_1$ meets the point $x_f$ and $R_q$ meets the point
$x_{f'}$  and such that $R_i\cap R_j$ is empty for $|i-j|>1$.
Then $\bd(e)=(\deg(E|_{R_1}),\dots,\deg(E|_{R_q}))$.
\end{enumerate}
\end{definition}

It is clear that for any chain-graph $(\tau,g,d,\bd)$ for rank $r$ 
there exists an $S_\tau$-labeled 
stable map to $\BGL_r$ whose combinatorial
type is $(\tau,g,d,\bd)$.

\section{The groupoid of stable maps to $\BGL_r$}
\label{gsm}

Let $S$ be a finite set and let $r\geq 1$. 
An 
{\em
$S$-labeled stable map over a scheme $T$ to $BGL_r$ 
}
is a tupel
$(C/T, x_i\ |\ i\in S, E)$, where $(C/T, x_i\ |\ i\in S)$ is a 
prestable $S$-pointed curve 
and $E$ is a vector bundle of rank $r$ on $C$ such that for each geometric
point $t$ of $T$  the fiber $(C_t, x_i(t)\ |\ i\in S, E_t)$ is a stable
map to $BGL_r$.

A morphism 
$$
(C/T, x_i\ |\ i\in S, E) \to 
(C'/T', x'_i\ |\ i\in S, E') 
$$
between two $S$-labeled stable maps to $\BGL_r$ consists of the following data:
\begin{enumerate}
\item
A morphism $T\to T'$,
\item
An isomorphism
$
f: C \isomto C'\times_{T'}T 
$
which for each $i\in S$  maps the section $x_i$ to the pull back 
of the section $x'_i$,
\item
An isomorphism
$
E \isomto f^*(\Oo_T\tensor_{\Oo_{T'}}E')
$.
\end{enumerate}

For integers $g\geq 0$ and $d$ 
let $\Mb_{g,S}(\BGL_r,d)$ be the groupoid (over the category of schemes)
whose objects are $S$-labeled stable maps 
$(C/T, x_i\ |\ i\in S, E)$
to $\BGL_r$ where $C/T$ is  a $T$-curve of genus $g$ and the restriction
of the vector bundle $E$ to each fiber of $C\to T$ is of degree $d$.

\section{Generalized isomorphisms and Gieseker vector bundle data}
\label{gi}

This section and the following are preparatory to
section \ref{bm} below, where we will define an analogue of the groupoid 
$\Mb(V,\tau)$ and the boundary morphism $\Mb(V,\tau)\to\Mb_{g,S}(V,\beta)$
(cf. \cite{Manin} V. §4.7) for the case when the variety $V$ is
replaced by the stack $\BGL_r$. 

Let $F_1$ and $F_2$ be two $r$-dimensional $k$-vector spaces.
In \cite{kgl} we have defined a smooth compactification $\KGL(F_1,F_2)$ of
the variety $\Isom(F_1,F_2)\isomorph\GL_r$ of linear isomorphisms from
$F_1$ to $F_2$
which has properties analogous to the so called 
``wonderful compactification'' of adjoint linear groups
introduced by De Concini and Procesi. 
We call the $k$-valued points of $\KGL(F_1,F_2)$
{\em generalized isomorphisms} from $F_1$ to $F_2$.
The complement 
of $\Isom(F_1,F_2)$ inside $\KGL(F_1,F_2)$ is a normal crossing divisor 
with smooth irreducible components 
$$
Y_i=Y_i(F_1,F_2)
\quad\text{and}\quad
Z_i=Z_i(F_1,F_2)
\quad (i\in[0,r-1])
\quad.
$$
For subsets $I,J$ of the set $[0,r-1]$ we write $\Ob_{I,J}(F_1,F_2)$ for
the intersection of all $Z_i$ and $Y_j$, where $i$ and $j$ run through the
set $I$ and $J$ respectively.
The variety $\Ob_{I,J}(F_1,F_2)$ is non-empty if and only if 
$$
\min(I)+\min(J)\geq r
$$ 
(here we use the convention that the minimum of the empty subset 
of $[0,r-1]$ is $r$).
A pair $(I,J)$ of subsets of $[0,r-1]$ which
satisfies this property is called a {\em GI-type} (for the rank $r$).
A $k$-rational point of $\Ob_{I,J}(F_1,F_2)$ is called a 
{\em generalized isomorphism of type $(I,J)$}.
There is an isomorphism 
$$
\inv:
\left\{
\begin{array}{ccc}
\KGL(F_1,F_2) &\to& \KGL(F_2,F_1) \\
\Phi &\mapsto& \Phi^{-1}
\end{array}
\right.
$$
with $\inv^2=\id$,
whose restriction to $\Isom(F_1,F_2)$ is the inversion.
We have $\inv(Y_i(F_1,F_2))=Z_i(F_2,F_1)$ for all $i\in[0,r-1]$.
Thus if $\Phi$ is a generalized isomorphism of type $(I,J)$
then $\Phi^{-1}$ is a generalized isomorphism of type $(J,I)$.

\begin{definition}
Let $(D,x_1,x_2)$ be a two-pointed prestable curve over a field $k$.
A {\em Gieseker vector bundle datum} on $(D,x_1,x_2)$ is a pair
$(F,\Phi)$ where $F$ is a vector bundle on $D$ with fibers $F_1$ and $F_2$
in $x_1$ and $x_2$ respectively,
and $\Phi$ is a generalized isomorphism from $F_1$ to $F_2$.
\end{definition}

Let $(D,x_1,x_2)$ be a two-pointed prestable curve and let 
$C'$ be the prestable curve which arises by identifying the
two points $x_1$ and $x_2$. Thus we have a morphism $D\to C'$
which is the partial normalization of the curve $C'$ at a singular
point $x$.

In our paper \cite{degeneration} we have shown that 
a Gieseker vector bundle datum $(F,\Phi)$ on $D$ induces a Gieseker vector
bundle $(C\to C',E)$ on $C'$.
Moreover the construction of $(C\to C',E)$ yields a distinguished
singular point $y\in C$ which is mapped to the point $x$,
and one can recover the data 
$(D,x_1,x_2,F,\Phi)$ from the data $(C\to C',E,y)$.

Let $(C\to C',E)$ be the Gieseker vector bundle associated to 
the Gieseker vector bundle datum $(F,\Phi)$.
Let $R$ be the fiber over $x$ of the morphism $C\to C'$.
Then $R$ is either a single point or a chain of projective lines. 
In the second case we may consider $D$ as a subscheme of $C$
and may number the components $R_1,\dots,R_m$ of $R$ 
such that $R_1$ meets $D$ in $x_1$ and $R_m$ meets $D$ in $x_2$
and $R_i\cap R_j=\emptyset$ if $|i-j|>1$.
Let $\bd=\emptyset$ if $R$ is reduced to a point and 
$\bd=(\deg(E|_{R_1}),\dots,\deg(E|_{R_m}))$ else.

Now what can be said about the chain-type $\bd$,
if the generalized isomorphism $\Phi$ is chosen in $\Ob_{I,J}(F_1,F_2)$
for some GI-type $(I,J)$?
To answer this question we define a mapping
$$
\{ \text{GI-types}\} \to \{\text{chain-types}\}
$$
as follows: Let $(I,J)$ be a GI-type, with 
$I=\{ i_1,\dots,i_p\}$,
$J=\{ j_1,\dots,j_q\}$
where $p,q\geq 0$ and  
\begin{eqnarray*}
&0\leq i_1\leq \dots \leq i_p<i_{p+1}:= r& \\
&0\leq j_1\leq \dots \leq j_q<j_{q+1}:= r&
\quad.
\end{eqnarray*}
The GI-type $(I,J)$ is mapped to the chain-type 
$(d_1,\dots,d_{p+q})$, where
$$
d_m=
\left\{
\begin{array}{ll}
i_{m+1}-i_m   &  \text{if $m\in[1,p]$} \\
j_{p+q-m+2}-j_{p+q-m+1} & \text{if $m\in[p+1,p+q]$}
\end{array}
\right.
$$
The mapping can be visualized by the following picture
$$
\overunderbraces%
{}
{
\overset{0}{|}
&
\dots\dots
&
\dots
&
\overset{i_1}{|}
&
\dots
&
\overset{i_2}{|}
\dots\dots
\overset{i_p}{|}
&
\dots
&
\overset{r}{|}
\overset{r}{|}
&
\dots
&
\overset{j_q}{|}
\dots\dots
\overset{j_2}{|}
&
\dots
&
\overset{j_1}{|}
&
\dots
&
\dots\dots
&
\overset{0}{|}
}{ 
& &
& &
\br{1}{d_1}
& &
\br{1}{d_p}
& &
\br{1}{d_{p+1}}
& &
\br{1}{d_{p+q}}
& &
& &
}
$$

The answer to the above question can now be formulated as follows:
If $\Phi$ is generic in $\Ob_{I,J}(F_1,F_2)$
(more precisely, if $\Phi$ is contained in $\Ob_{I,J}(F_1,F_2)$ but not in
$\Ob_{I',J'}(F_1,F_2)$ for $I'$ strictly larger than $I$ 
or $J'$ strictly larger than $J$) 
then the chain-type $\bd$ is the image of the GI-type $(I,J)$ by the mapping
$$
\{ \text{GI-types}\} \to \{\text{chain-types}\}
$$
defined above.

\section{The groupoid of Gieseker vector bundle data}
\label{DB}

Let $T$ be a scheme and let $F_1$ and $F_2$ be two
vector bundles of rank $r$ on $T$. 
In \cite{kgl} we have defined 
a smooth $T$-scheme $\KGL(F_1,F_2)$, whose fiber over at a point $t\in T$ is a 
compactification of the variety $\Isom(F_{1,t},F_{2,t})$, which is isomorphic
to $\KGL(F_{1,t},F_{2,t})$.
A section of $\KGL(F_1,F_2)\to T$ is called a {\em generalized isomorphism 
from $F_1$ to $F_2$}.
For each GI-type $(I,J)$ there exists a natural closed subscheme 
$\Ob_{I,J}(F_1,F_2)$ 
of $\KGL(F_1,F_2)$ whose fiber over $t\in T$ is the variety
$\Ob_{I,J}(F_{1,t},F_{2,t})$.
A section of $\Ob_{I,J}(F_1,F_2)\to T$ is called a 
{\em generalized isomorphism of type $(I,J)$ from $F_1$ to $F_2$}.
Also the isomorphism 
$\inv:\KGL(F_1,F_2)\isomto\KGL(F_2,F_1)$, $\Phi\mapsto\Phi^{-1}$ 
generalizes to the relative situation and has the properties 
$\inv^2=\id$ and $\inv(\Ob_{I,J}(F_1,F_2))=\inv(\Ob_{J,I}(F_2,F_1)$.
All these objects commute (in an obvious sense) with base change 
$T'\to T$.

\begin{definition}
Let $(D/T,x_1,x_2)$ be a two-pointed prestable curve over a scheme $T$
and let $(I,J)$ be a GI-type for a rank $r\in\N$.
A {\em Gieseker vector bundle datum of rank $r$ and type 
   $(I,J)$ on $(D/T,x_1,x_2)$}
is a pair $(F,\Phi)$, where $F$ is a vector bundle of rank $r$ on 
$D$ and $\Phi$ is a section of $\Ob_{I,J}(x_1^*F,x_2^*F)\to T$.
\end{definition}

An isomorphism 
$$
(F,\Phi)\isomto(F',\Phi')
$$
between two Gieseker vector bundle data of type $(I,J)$ on 
$(D/T,x_1,x_2)$
is an isomorphism $F\isomto F'$ such
that the induced isomorphism
$\Ob_{I,J}(x_1^*F,x_2^*F)\isomto\Ob_{I,J}(x_1^*F',x_2^*F')$
carries $\Phi$ to $\Phi'$.

We denote by $\GVBD_{I,J}$ the groupoid whose objects are Gieseker vector
bundle data of rank $r$ and type $(I,J)$ on $(D/T,x_1,x_2)$ and whose 
arrows are isomorphisms between Gieseker vector bundle data.

Let $(I,J)$ and $(I',J')$ be two GI-types which are mapped to the
same chain-type. 
Then by similar arguments as used in the proof of 
\cite{factorization} 4.2 it can be shown that 
there is a canonical isomorphism of groupoids
$$
\beta^{I,J}_{I',J'}:
\GVBD_{I,J}\isomto\GVBD_{I',J'}
\quad.
$$
Moreover, we have $\beta^{I,J}_{I,J}=\id$ and if 
$(I'',J'')$ is a third GI-type with the same associated chain-type
as $(I,J)$ and $(I',J')$, then 
$\beta^{I,J}_{I'',J''}=\beta^{I',J'}_{I'',J''}\comp\beta^{I,J}_{I',J'}$.

Let $\bd$ be a chain-type.
The family of isomorphisms 
$(\beta^{I,J}_{I',J'})_{(I,J), (I',J')\mapsto\bd}$ 
enables us to define a new groupoid $\GVBD_{\bd}$ as follows.
The objects of $\GVBD_{\bd}$ are the union of all objects of the 
groupoids  $\GVBD_{I,J}$ where $(I,J)$ runs through all GI-types which are
mapped to the chain-type $\bd$.
A morphism from an object $(F,\Phi)$ of $\GVBD_{I,J}$ to an object
$(F',\Phi')$ of $\GVBD_{I',J'}$ is a family $(f_{A,B})_{(A,B)\mapsto\bd}$
where 
$$
f_{A,B}:\beta^{I,J}_{A,B}(F,\Phi)\to\beta^{I',J'}_{A,B}(F',\Phi')
$$
is a morphism in $\GVBD_{A,B}$ such that 
$\beta^{A,B}_{A',B'}(f_{A,B})=f_{A',B'}$
for all GI-types $(A,B)$, $(A',B')$ mapped to $\bd$.

\begin{remark}
Let $(I,J)$ be a GI-type and let $\bd$ be its associated 
chain-type. Then the inclusion $\GVBD_{I,J}\injto \GVBD_{\bd}$
is an equivalence of categories.
The advantage of $\GVBD_{\bd}$ over $\GVBD_{I,J}$ is that its
definition depends only on the chain-type $\bd$ and 
does not involve the choice of of a GI-type over $\bd$.
\end{remark}

\begin{definition}
An object $(F,\Phi)$ of the groupoid $\GVBD_{\bd}$ is called
a {\em Gieseker vector bundle datum of type $\bd$} on the two-pointed curve 
$(D/T,x_1,x_2)$.
\end{definition}

Let $C'$ be the prestable $T$-curve constructed from $D$ by
identifying the two sections $x_1$ and $x_2$ and let $x$ be the
section of $C'$ given by the composition 
$T\overset{x_i}{\to} D\to C'$.
Let $\GVB$ be the following groupoid:
The objects of $\GVB$ are Gieseker vector bundles $(C\to C',E)$ 
of rank $r$ on $C'$ such that $C\to C'$ is an isomorphism over
the complement of $x(T)$.
An arrow $(C_1\to C',E_1)\to(C_2\to C',E_2)$ between two such
objects consists of a $C'$-isomorphism
$f:C_1\isomto C_2$ and an isomorphism $E_1\isomto f^*E_2$.

In section \ref{gi} we have already mentioned that a Gieseker vector
bundle datum on a two-pointed prestable curve over a field induces a
Gieseker vector bundle on the prestable curve obtained by identifying
the two marked points. This construction works also in families over
an arbitrary base scheme $T$ (cf. \cite{degeneration}) and gives for
each GI-type $(I,J)$ a natural morphism of groupoids
$$
\nu_{I,J}: \GVBD_{I,J}\to\GVB
\quad.
$$

The morphisms $\nu_{I,J}$ are compatible with the isomorphisms 
$\beta^{I,J}_{I',J'}$ in the sense that all triangles
$$
\xymatrix{
\GVBD_{I,J} 
\ar[dr]_{\text{$\nu_{I,J}$}} 
\ar[rr]^{\text{$\beta^{I,J}_{I',J'}$}} 
& &
\GVBD_{I',J'} \ar[dl]^{\text{$\nu_{I',J'}$}}  \\
& \GVB
}
$$
commute. 

It follows that for each chain-type $\bd$ the family 
$(\nu_{I,J})_{(I,J)\mapsto\bd}$
induces a morphism of groupoids 
$$
\nu_{\bd}:\GVBD_{\bd}\to \GVB
\quad.
$$

\section{Boundary morphisms}
\label{bm}

In this section we will define for
each chain-graph (cf. section \ref{ct})
$$
\btau=
(\tau,\ g_v\ |\ v\in V_\tau,\ d_v\ |\ v\in V_\tau,\ 
      \bd_e\ |\ e\in E_\tau^{\ori})
$$ 
a groupoid (over the category
of schemes) $\Mb(\BGL_r,\btau)$ together with a boundary morphism
$$
\nu_{\btau}: \Mb(\BGL_r,\btau)\to\Mb_{g,S}(\BGL_r,d)
\quad,
$$
where $g=\sum g_v + \dim H_1(||\tau||)$ is the genus of the 
combinatorial type, $d=\sum d_v + \sum |\bd_e|$ is its total degree
and $S=S_\tau$ is the set of its tails.

We start with an auxiliary definition.
\begin{definition}
A {\em GI-graph} (for rank $r$) is a tupel
$$
\bgamma=
(\tau, \ g_v\ |\ v\in V_\tau,\ \delta_v\ |\ v\in V_\tau,\
       I_f\ |\ f\in F_\tau\setminus S_\tau) 
$$
where $(\tau, g_v\ | v\in V_\tau)$ is a stable modular graph,
where the $\delta_v$ are integers and where 
for each oriented edge $(f,f')$ the pair $(I_f,I_{f'})$ is
a GI-type for the rank $r$.
\end{definition}

Before defining the groupoid $\Mb(\BGL_r,\btau)$ we will first define
the groupoid $\Mb(\BGL_r,\bgamma)$ for a GI-graph
$\bgamma$.
We will see below that to $\bgamma$ there is associated a 
chain-graph $\btau$, and 
the relationship between the groupoids 
$\Mb(\BGL_r,\bgamma)$ and 
$\Mb(\BGL_r,\btau)$
is analogous to the relationship between 
$\GVBD_{I,J}$ and $\GVBD_{\bd}$ described in section \ref{DB}.

Let 
$$
\bgamma=
(\tau, \ g_v\ |\ v\in V_\tau,\ \delta_v\ |\ v\in V_\tau,\
       I_f\ |\ f\in F_\tau\setminus S_\tau) 
$$
be a GI-graph and let $T$ be a scheme. 
An object of $\Mb(\BGL_r,\bgamma)$ over $T$ is a family
$$
((C_v/T ,\ (x_i)_{i\in F_\tau(v)},\ \F_v)\ |\ v\in V_\tau
\ \ ;\ \
 \Phi_e\ |\ e\in E_\tau^{\ori})
$$
where
\begin{enumerate}
\item
$(C_v/T ,\ (x_i)_{i\in F_\tau(v)},\ \F_v)$
is an $F_\tau(v)$-labeled stable map over $T$ to $\BGL_r$ such
that $C_v/T$ is a curve of genus $g$ over $T$.
\item
For every oriented edge $e=(f,f')$ of $\tau$ 
with $f\in F_\tau(v)$ and $f'\in F_\tau(v')$,
$\Phi_e$ is a generalized isomorphism of type $(I_f,I_{f'})$ from 
$x_f^*\F_v$ to $x_{f'}^*\F_{v'}$.
\item
For every oriented edge $e$ 
the equality $\Phi_{j(e)}=\Phi_e^{-1}$ holds.
\item
For each vertex $v\in V_\tau$ we have
$
\deg{\F_v}=\delta_v
$.
\end{enumerate}

An arrow 
\begin{eqnarray*}
&
((C_v/T ,\ (x_i)_{i\in F_\tau(v)},\ \F_v)\ |\ v\in V_\tau
\ \ ;\ \
 \Phi_e\ |\ e\in E_\tau^{\ori})
&
\\
\to
&
((C'_v/T ,\ (x'_i)_{i\in F_\tau(v)},\ \F'_v)\ |\ v\in V_\tau
\ \ ;\ \
 \Phi'_e\ |\ e\in E_\tau^{\ori})
&
\end{eqnarray*}
between two objects of $\Mb(\BGL_r,\bgamma)$ over $T$ is
a family of isomorphisms
$$
(C_v/T ,\ (x_i)_{i\in F_\tau(v)},\ \F_v)
\isomto
(C'_v/T ,\ (x'_i)_{i\in F_\tau(v)},\ \F'_v)
$$
of stable maps over $T$ to $\BGL_r$ (cf. section \ref{gsm}), where
$v$ runs through the set $V_\tau$ of vertices of $\tau$,
such that for each oriented edge $e=(f_1,f_2)$ 
with $f_1\in F_\tau(v_1)$, $f_2\in F_\tau(v_2)$
the induced isomorphism
$$
\KGL(x_{f_1}^*\F_{v_1}, x_{f_2}^*\F_{v_2}) \isomto
\KGL((x'_{f_1})^*\F'_{v_1}, (x'_{f_2})^*\F'_{v_2})
$$
maps the generalized isomorphism $\Phi_e$ to $\Phi'_e$.

If $T'\to T$ is a morphism of schemes then an object of
$\Mb(\BGL_r,\bgamma)$ over $T$ pulls back to one over $T'$ 
and an arrow between objects over $T$ induces one between
their respective pull back over $T'$.
This defines $\Mb(\BGL_r,\bgamma)$ 
as a groupoid over the category of schemes.

To a GI-graph $\bgamma$ we associate the chain-graph
$$
\btau=
(\tau,\ g_v\ |\ v\in V_\tau,\ d_v\ |\ v\in V_\tau,\ 
      \bd_e\ |\ e\in E_\tau^{\ori})
$$ 
where 
$
d_v=\delta_v-\sum_{f\in F_\tau(v)}(r-\min(I_f))
$,
and for each oriented edge $e=(f,f')$ the chain-type $\bd_e$ is the image 
of the GI-type $(I_f,I_{f'})$.
This yields a mapping
$$
\left\{
\begin{array}{l}
\text{GI-graphs}
\end{array}
\right\}
\To
\left\{
\begin{array}{l}
\text{chain-graphs}
\end{array}
\right\}
$$
which obviously is surjective and finite to one.

Let $\btau$ be a chain-graph and let $\bgamma$ and $\bgamma'$ be two 
GI-graphs mapped to $\btau$.
Using arguments similar as in the proof of \cite{factorization} 4.2,
one shows that there is a canonical isomorphism  
$$
\beta^{\bgamma}_{\bgamma'}: \Mb(\BGL_r,\bgamma)\to\Mb(\BGL_r,\bgamma')
\quad.
$$
We have $\beta^{\bgamma}_{\bgamma}=\id$ and 
$
\beta^{\bgamma'}_{\bgamma''}\comp\beta^{\bgamma}_{\bgamma'}=
\beta^{\bgamma}_{\bgamma''}
$
for any further GI-graph $\bgamma''$ mapped to $\btau$.

We use 
the family $(\beta^{\bgamma}_{\bgamma'})_{\bgamma,\bgamma'\mapsto\btau}$
to define the groupoid $\Mb(\BGL_r,\btau)$ over the category of schemes.
The construction is completely analogous to the construction
of $\GVBD_{\bd}$ in section \ref{DB}:
The objects of $\Mb(\BGL_r,\btau)$ are the union of objects
of $\Mb(\BGL_r,\bgamma)$ where $\bgamma$ runs through all
GI-graphs mapped to $\btau$.
Let $\E$ and $\E'$ be two objects in 
$\Mb(\BGL_r,\btau)$ and assume that 
$\E$ belongs to $\Mb(\BGL_r,\bgamma)$ and
$\E'$ to $\Mb(\BGL_r,\bgamma')$.
An arrow of $\Mb(\BGL_r,\btau)$ from 
$\E$ to $\E'$ is then a family of arrows
$$
f_{\bpi}: \beta^{\bgamma}_{\bpi}(\E)\to\beta^{\bgamma'}_{\bpi}(\E')
$$
in $\Mb(\BGL_r,\bpi)$
where $\bpi$ runs through all GI-graphs over $\btau$
such that for all pairs $\bpi$, $\bpi'$ of GI-graphs
over $\btau$ we have 
$\beta^{\bpi}_{\bpi'}(f_{\bpi})=f_{\bpi'}$.

\begin{remark}
Of course for each GI-graph $\bgamma$ mapped to a chain-graph
$\btau$ the inclusion 
$\Mb(\BGL_r,\bgamma)\injto\Mb(\BGL_r,\btau)$
is an equivalence of categories
and the reason why we defined
$\Mb(\BGL_r,\btau)$ in the way we did
is that we want to make apparent that it does not depend
on any choice of a GI-graph $\bgamma$ over $\btau$.
\end{remark}

The boundary morphism
$$
\nu_{\btau}: \Mb(\BGL_r,\btau)\to \Mb_{g,S}(\BGL_r, d)
$$
is now defined as follows:
With the help of the constructions in \cite{degeneration}
we get a canonical morphism 
$\nu_{\bgamma}: \Mb(\BGL_r,\bgamma)\to\Mb_{g,S}(\BGL_r, d)$
for every GI-graph $\bgamma$ over $\btau$ such that
for any two GI-graphs $\bgamma$ and $\bgamma'$ the diagram
$$
\xymatrix{
\text{$\Mb(\BGL_r,\bgamma)$}
\ar[rd]_{\text{$\nu_{\bgamma}$}}
\ar[rr]^{\text{$\beta^{\bgamma}_{\bgamma'}$}}
& &
\text{$\Mb(\BGL_r,\bgamma')$}
\ar[ld]^{\text{$\nu_{\bgamma'}$}}
\\
& \text{$\Mb_{g,S}(\BGL_r, d)$} &
}
$$
Therefore the family $(\nu_{\bgamma})_{\bgamma\mapsto\btau}$
induces a morphism 
$
\nu_{\btau}: \Mb(\BGL_r,\btau)\to \Mb_{g,S}(\BGL_r, d)
$.


\begin{thebibliography}{99}
\bibitem[AV]{AV}
         D. Abramovich, A. Vistoli:
         Compactifying the space of stable maps.
         J. Amer. Math. Soc.  15  (2002),  no. 1, 27--75
\bibitem[ACV]{ACV}
         D. Abramovich, A. Corti, A. Vistoli:
         Twisted bundles and admissible covers.
         Special issue in honor of Steven L. Kleiman.
         Comm. Algebra 31 (2003), no. 8, 3547--3618.
\bibitem[G]{Gieseker}
         Gieseker:
         A degeneration of the moduli space of stable bundles.
         J. Differential Geom. 19 (1984), no. 1, 173-206.
\bibitem[K1]{kgl}
         I. Kausz:
         A Modular Compactification 
         of the General Linear Group,
         Documenta Math. 5 (2000) 553-594
\bibitem[K2]{degeneration}
         I. Kausz:
         A Gieseker Type Degeneration of Moduli Stacks of Vector
         Bundles on Curves,
         math.AG/0201197
\bibitem[K3]{factorization}
         I. Kausz:
         A canonical decomposition of generalized theta functions on the
         moduli stack of Gieseker vector bundles.
         math.AG/0305034
\bibitem[K4]{ramified_bundles}
         I. Kausz:
	 Twisted vector bundles on pointed nodal curves.  
	 Proc. Indian Acad. Sci. Math. Sci.  115  (2005),  no. 2, 147--165.
\bibitem[Li]{Li}
         J. Li:
	 Moduli spaces associated to a singular variety and the 
	 moduli of bundles over universal curves.
	 Vector bundles and representation theory (Columbia, MO, 2002), 57--74,
	 Contemp. Math., 322,
	 Amer. Math. Soc., Providence, RI, 2003. 
\bibitem[M]{Manin}
         Y. I. Manin:
         Frobenius Manifolds, Quantum Cohomology, and Moduli Spaces.
         AMS Colloquium Publications (1999).
\bibitem[Sch]{Sch}
         A. Schmitt:
	 The Hilbert compactification of the universal moduli space 
	 of semistable vector bundles over smooth curves.  
	 J. Differential Geom. 66 (2004), no. 2, 169--209.
\end{thebibliography}
\end{document}